\documentclass[11pt,reqno]{amsart}

\usepackage{amsmath,amssymb,amsthm,mathtools}
\usepackage{booktabs,array}
\usepackage[hidelinks]{hyperref}

\newtheorem{theorem}{Theorem}[section]
\newtheorem{proposition}[theorem]{Proposition}
\newtheorem{lemma}[theorem]{Lemma}

\theoremstyle{definition}
\newtheorem{definition}[theorem]{Definition}
\newtheorem{example}[theorem]{Example}

\theoremstyle{remark}
\newtheorem{remark}[theorem]{Remark}

\DeclareMathOperator{\lcm}{lcm}
\DeclareMathOperator{\Par}{Par}

\title[Prime Factorization and Exact Covering Systems]{Prime Factorization and the Structure of Exact Covering Systems}

	\author{Snir Kadosh}
\address{Department of Mathematics, ORT Braude College, 2161002 Karmiel, Israel}
\email{Snir.Kadosh@e.braude.ac.il}

\author{Ofir Schnabel}
\address{Department of Mathematics, ORT Braude College, 2161002 Karmiel, Israel}
\email{ofirsch@braude.ac.il}

\begin{document}
	
	\begin{abstract}
		An exact covering system (ECS) is a finite collection of arithmetic progressions that partitions $\mathbb{Z}$. We study two structural questions governed by the prime factorization of the system's length $L = \lcm(n_1,\dots,n_k)$. First, using generating functions and an asymmetric minimal vanishing sum of roots of unity due to Lam and Leung, we construct an explicit infinite family $C_{r,p,q}$ of irreducible ECSs of length $rpq$ for every triple of distinct primes $r, p, q$, generalizing a classical example of length $30$. We also exhibit an irreducible ECS of length $30$ outside this family, showing that $C_{r,p,q}$ does not exhaust all irreducible systems. Second, we completely classify the lengths $n > 1$ for which every multiset of moduli satisfying the three classical necessary conditions--the density condition, the pairwise non-coprime condition, and the Chen--Porubsk\'y multiplicity condition--is realizable as an ECS. Specifically, we prove that all such candidate multisets of length $n > 1$ are realizable if and only if $n$ has the form $p^k$, $pq$, or $p^2q$ for distinct primes $p \neq q$ and $k \geq 1$.
	\end{abstract}
	
	\maketitle
	
	\section{Introduction}
	
	For integers $a$ and $n \geq 1$, denote the arithmetic progression (or congruence class) with residue $a$ and modulus $n$ by
	\[
	a(n) = a + n\mathbb{Z} = \{a + jn : j \in \mathbb{Z}\},
	\]
	taking canonical representatives $0 \leq a < n$. A finite collection $\mathcal{A} = \{a_i(n_i)\}_{i=1}^k$ is a \emph{covering system} if $\bigcup_{i=1}^k a_i(n_i) = \mathbb{Z}$, and an \emph{exact covering system} (ECS), or \emph{disjoint covering system}, if every integer belongs to a unique class in $\mathcal{A}$. Introduced by Erd\H{o}s~\cite{Erdos1950} in connection with integers of the form $2^k + p$, covering systems have generated an extensive literature; see~\cite{Porubsky1981,PorubskySchonheim2002} for comprehensive surveys.
	
	The \emph{length} of $\mathcal{A} = \{a_i(n_i)\}_{i=1}^k$ is $L(\mathcal{A}) = \lcm(n_1, \dots, n_k)$. Since membership in each class depends only on residues modulo $L(\mathcal{A})$, exactness is equivalent to partitioning the interval $[0, L(\mathcal{A}) - 1]$. Counting covered elements in $[0, L(\mathcal{A}) - 1]$ yields the \emph{density condition}
	\begin{equation}\label{eq:density}
		D(\mathcal{A}) := \sum_{i=1}^k \frac{1}{n_i} = 1.
	\end{equation}
	Moreover, pairwise disjointness implies $\gcd(n_i, n_j) > 1$ for all $i \neq j$ by the Chinese Remainder Theorem~\cite{Korec1984}. A classical theorem of Mirsky, Newman, Davenport, and Rad\'o~\cite{Erdos1952} shows that the largest modulus of any nontrivial ECS occurs at least twice; more generally, Chen and Porubsk\'y~\cite{ChenPorubsky1995} proved that if $p_1 < \cdots < p_m$ are the distinct prime divisors of the largest modulus $n_k$, then its multiplicity $t$ satisfies
	\begin{equation}\label{eq:chen_porubsky}
		t = t_1 p_1 + \cdots + t_m p_m \qquad \text{for some } t_1, \dots, t_m \in \mathbb{Z}_{\geq 0}.
	\end{equation}
	
	A fundamental operation on ECSs is \emph{prime refinement}: for a prime $p$, replacing a single class $a(n)$ by the $p$ classes $\{(a + jn)(pn)\}_{j=0}^{p-1}$ preserves exactness. An ECS is \emph{reducible} if it can be obtained from a coarser ECS by a nonempty finite sequence of prime refinements, and \emph{irreducible} otherwise. (Equivalently, an ECS is irreducible if no $p$ classes of a common modulus $pn$ share a residue modulo $n$.) An ECS is \emph{natural} (NECS) if it arises from the trivial system $\{0(1)\}$ via a finite sequence of prime refinements~\cite{Goulden2018}. Korec~\cite{Korec1984} and Schnabel~\cite{Schnabel2015} proved that every nontrivial ECS whose length has at most two distinct prime divisors is natural, and hence reducible. Thus any nontrivial irreducible ECS must have a length divisible by at least three distinct primes.
	
	In this paper, we establish two main theorems illustrating how the prime factorization of the length governs both irreducibility and realizability. First, translating exactness into vanishing sums of roots of unity via generating functions and invoking an asymmetric minimal vanishing sum of Lam and Leung~\cite{LamLeung2000}, we construct an infinite family of irreducible ECSs parametrized by arbitrary triples of distinct primes.
	
	\medskip
	\noindent\textbf{Theorem A.}
	\emph{For any three distinct primes $r, p, q$, the collection
		\begin{align*}
			C_{r,p,q} &= \bigl\{ rk(rp) : 1 \leq k \leq p-1 \bigr\} \cup \{0(pq)\} \cup \bigl\{ k(rq) : 1 \leq k < rq,\; r \nmid k,\; q \nmid k \bigr\} \\
			&\quad \cup \bigl\{ rpk(rpq) : 1 \leq k \leq q-1 \bigr\}\\
			& \cup \bigl\{ bq(rpq) : 1 \leq k \leq p-1,\; 1 \leq m \leq r-1,\; b \equiv rk + pm \pmod{rp} \bigr\}
		\end{align*}
		is an irreducible exact covering system of length $rpq$.}
	\medskip
	
	For $(r,p,q)=(2,3,5)$, Theorem~A recovers the classical length-$30$ irreducible system recorded in~\cite{Schnabel2015}. To the best of our knowledge, apart from the family \(C_{r,p,q}\) and its expansions, all currently known irreducible ECSs appear to be isolated examples. We also exhibit an irreducible ECS of length $30$ lying outside the family $C_{r,p,q}$. To the best of our knowledge, this is the first explicit irreducible ECS whose irreducibility is not explained by the Lam--Leung minimal vanishing sum construction, demonstrating that the class of irreducible ECSs is substantially richer than the family provided by Theorem~A. 
	
	Whereas Theorem~A concerns the existence of irreducible exact covering systems, Theorem~B addresses the complementary question of which multisets of moduli can occur in an exact covering system.
	
	We study the realizability problem for moduli multisets: given a multiset $M = \{n_1, \dots, n_k\}$ of positive integers, when do there exist residues $a_1, \dots, a_k$ such that $\{a_i(n_i)\}_{i=1}^k$ is an ECS? For each $n\ge1$, let $\mathcal C(n)$ denote the set of candidate multisets of divisors of $n$ with least common multiple $n$ satisfying the density condition~\eqref{eq:density}, the pairwise non-coprime condition, and~\eqref{eq:chen_porubsky}. Let $\mathcal R(n)\subseteq\mathcal C(n)$ denote the subset of realizable multisets. The multiset \[ M=\{2,4,6,24,24\} \] belongs to $\mathcal C(24)$ but is not realizable, showing that these three necessary conditions are not sufficient in general. We completely classify the lengths $n>1$ for which the classical necessary conditions are also sufficient for realizability.
		
	\noindent\textbf{Theorem B.} \emph{ For every integer $n>1$, \[ \mathcal R(n)=\mathcal C(n) \] if and only if \[ n=p^k,\qquad n=pq,\qquad\text{or}\qquad n=p^2q, \] where $p\neq q$ are primes and $k\ge1$. }
	
	\medskip 
Section~\ref{sec:preliminaries} reviews basic operations on exact covering systems and develops the generating-function framework used throughout the paper. Section~\ref{sec:irreducible} proves Theorem~A and discusses irreducible ECSs beyond the Lam--Leung family. Section~\ref{sec:classification} establishes Theorem~B.
	
	\section{Preliminaries and Generating Functions}\label{sec:preliminaries}
	
	We record two basic operations on covering systems~\cite{Goulden2018} and a criterion for disjointness.
	
	\begin{definition}\label{def:operations} Let $\mathcal A=\{a_i(n_i)\}_{i=1}^k$ be a covering system. For $r\in\mathbb Z$ and $s\ge1$, define \[ s\mathcal A+r = \{(r+sa_i)(sn_i)\}_{i=1}^k. \] We call $s\mathcal A+r$ the expansion of $\mathcal A$ by $r$ and $s$. \end{definition}
	
	The operation preserves exactness. When $s=1$ it merely translates all residues by $r$, while for $s\ge2$ it produces an exact partition of the congruence class $r(s)$ into classes whose moduli are scaled by $s$.
	
	\begin{lemma}\label{lem:intersection}
		Two arithmetic progressions $a(m)$ and $b(n)$ satisfy $a(m) \cap b(n) \neq \emptyset$ if and only if $\gcd(m, n) \mid (a - b)$.
	\end{lemma}
	
	\begin{proof}
		The intersection $a(m) \cap b(n)$ is nonempty if and only if the system $x \equiv a \pmod{m}$, $x \equiv b \pmod{n}$ is solvable, which by the generalized Chinese Remainder Theorem holds if and only if $a \equiv b \pmod{\gcd(m, n)}$.
	\end{proof}
	
	To study prime refinements algebraically, we associate to an ECS $\mathcal{A} = \{a_i(n_i)\}_{i=1}^k$ (with $0 \leq a_i < n_i$) its generating-function identity. Since $\mathcal{A}$ partitions $\mathbb{Z}_{\geq 0}$, summing $z^m$ over $m \geq 0$ yields
	\[
	\sum_{i=1}^k \frac{z^{a_i}}{1 - z^{n_i}} = \frac{1}{1 - z}, \qquad |z| < 1.
	\]
	Setting $L = L(\mathcal{A})$ and multiplying both sides by $1 - z^L$ gives the polynomial identity
	\[
	\sum_{i=1}^k z^{a_i} \sum_{j=0}^{L/n_i - 1} z^{j n_i} = \sum_{m=0}^{L-1} z^m.
	\]
	Let $\zeta = e^{2\pi i / L}$ be a primitive $L$-th root of unity. Evaluating at $z = \zeta$, the right-hand side vanishes, while for every $n_i < L$ the factor $\sum_{j=0}^{L/n_i - 1} (\zeta^{n_i})^j$ vanishes because $\zeta^{n_i}$ is a primitive $(L/n_i)$-th root of unity. If $L$ occurs as a modulus in $\mathcal{A}$, this isolates the classes of modulus $L$ (called \emph{singletons} relative to $[0, L-1]$):
	\begin{equation}\label{eq:vanishing_sum}
		\sum_{\substack{1 \leq i \leq k \\ n_i = L}} \zeta^{a_i} = 0.
	\end{equation}
	If $p\mid L$ is prime and $\mathcal A$ contains the complete $p$-refinement block \[ \{(a+jL/p)(L):0\le j\le p-1\}, \] then the corresponding roots of unity form the vanishing $p$-cycle \[ \sum_{j=0}^{p-1}\zeta^{a+jL/p} = \zeta^a\sum_{j=0}^{p-1}(\zeta^{L/p})^j = 0. \] Conversely, if~\eqref{eq:vanishing_sum} is a minimal vanishing sum, meaning that no nonempty proper subsum vanishes, then the classes of modulus $L$ cannot contain a complete $p$-cycle for any prime $p\mid L$. Consequently, no prime refinement among the classes of modulus $L$ can be reversed.
	\section{An Infinite Family of Irreducible ECSs}\label{sec:irreducible}
	Lam and Leung~\cite[Example~2.5]{LamLeung2000} constructed an asymmetric minimal vanishing sum of $rpq$-th roots of unity for any three distinct primes $r, p, q$. Choosing the classes of modulus $rpq$ so that their associated roots of unity realize this minimal vanishing sum, and then completing the remaining residue classes to obtain a covering, produces the family $C_{r,p,q}$.
	\begin{proof}[Proof of Theorem~A]
		Partition $C_{r,p,q}$ into five subcollections:
	\begin{align*} \mathcal{U} &= \{rk(rp):1\le k\le p-1\},\\ \mathcal{D} &= \{0(pq)\},\\ \mathcal{T} &= \{k(rq):1\le k<rq,\ r\nmid k,\ q\nmid k\},\\ \mathcal{S} &= \{rpk(rpq):1\le k\le q-1\},\\ \mathcal{R} &= \{bq(rpq):1\le k\le p-1,\ 1\le m\le r-1,\ 0\le b<rp,\ b\equiv rk+pm\pmod{rp}\}. \end{align*}
		\medskip
		\noindent\emph{Step 1: Coverage.}
		Split $\mathcal{D} = \{0(pq)\}$ into the $r$ classes $D_t = \{pqt(rpq)\}$ for $1 \leq t \leq r$. Taking the union of $\mathcal{S}$ with $D_r = \{0(rpq)\}$ gives
		\[ \mathcal S\cup D_r = \{\,rpk(rpq):0\le k\le q-1\,\} = 0(rp). \]
		Combining this with $\mathcal{U}$ yields
		\[ \mathcal U\cup\mathcal S\cup D_r = \{\,rk(rp):0\le k\le p-1\,\} = 0(r). \]
		Next, for each $1 \leq t \leq r-1$, let $\mathcal{R}_t \subseteq \mathcal{R}$ be the subset of classes corresponding to $m = t$. As $k$ runs over ${1, \dots, p-1}$, $rk + pt$ runs over all residues modulo $rp$ congruent to $pt \pmod{r}$ except $pt$ itself. Together with $D_t = \{pqt(rpq)\}$, we obtain
		\(\mathcal{R}_{t}\cup D_{t}=\{(rk+pt)q(rpq):0\le k\le p-1\}=tq(rq).\)
		Now fix $1 \leq t \leq r-1$ and take the union of ${tq(rq)}$ with the subset $\mathcal{T}_t = \{\,k(rq): k\equiv t \pmod r,\ q\nmid k\,\}.$ Every integer $x \in t(r)$ satisfies $x \equiv t \pmod{r}$ and $x \equiv b \pmod{q}$ for a unique $0 \leq b < q$. If $b = 0$, then $x \equiv tq \pmod{rq}$ by the Chinese Remainder Theorem, so $x \in tq(rq)$; if $1 \leq b \leq q-1$, then $x \equiv k \pmod{rq}$ for a unique $k$ with $k \equiv t \pmod{r}$ and $q \nmid k$, so $x \in \mathcal{T}_t$. Hence \[ \mathcal{R}_t \cup D_t \cup \mathcal{T}_t = t(r). \] Taking the union over $0 \leq t \leq r-1$ yields $\bigcup_{t=0}^{r-1} t(r) = \mathbb{Z}$, so $C_{r,p,q}$ is a covering system.
		\medskip
		\noindent\emph{Step 2: Exactness.}
		Counting the classes in each subcollection--namely $p-1$ of modulus $rp$, one of modulus $pq$, $\varphi(rq) = (r-1)(q-1)$ of modulus $rq$, and $(q-1) + (r-1)(p-1)$ of modulus $rpq$--we compute the density:
		\begin{align*} D(C_{r,p,q}) &= \frac{p-1}{rp} +\frac1{pq} +\frac{(r-1)(q-1)}{rq} +\frac{(r-1)(p-1)+q-1}{rpq} \\ &= \frac{(p-1)q+r+p(r-1)(q-1)+(r-1)(p-1)+q-1}{rpq} =1. \end{align*}
	Since $C_{r,p,q}$ is a covering system and its densities sum to $1$, the classes cannot overlap on a set of positive density. Hence $C_{r,p,q}$ is exact.
		\medskip
		\noindent\emph{Step 3: Irreducibility.}
		We show that no prime refinement can be reversed at any modulus occurring in $C_{r,p,q}$:
		\begin{itemize}
			\item \emph{Modulus $rp$:} Only $p-1$ classes of modulus $rp$ appear, so no $p$-refinement can be reversed. Every class in $\mathcal{U}$ has residue divisible by $r$; since $p \not\equiv 0 \pmod{r}$, no residue class modulo $p$ contains $r$ classes of modulus $rp$, ruling out a reverse $r$-refinement.
			\item \emph{Modulus $pq$:} Only the single class $0(pq)$ appears.
			\item \emph{Modulus $rq$:} A reverse $r$-refinement (resp.\ $q$-refinement) requires $r$ classes of modulus $rq$ sharing a residue modulo $q$ (resp.\ $q$ classes sharing a residue modulo $r$). By the Chinese Remainder Theorem, one of these $r$ (resp.\ $q$) classes has residue divisible by $r$ (resp.\ $q$), which is excluded from $\mathcal{T}$.
			\item \emph{Modulus $rpq$:} Let $\zeta = e^{2\pi i / (rpq)}$ and set $t = \zeta^{pq}$, $h = \zeta^{rq}$, and $s = \zeta^{rp}$, which are primitive roots of unity of orders $r$, $p$, and $q$, respectively. The singletons $\mathcal{S} \cup \mathcal{R}$ of modulus $rpq$ correspond to the vanishing sum
			\(V=\sum _{m=1}^{r-1}\sum _{k=1}^{p-1}t^{m}h^{k}+\sum _{a=1}^{q-1}s^{a}=\left(\sum _{m=1}^{r-1}t^{m}\right)\left(\sum _{k=1}^{p-1}h^{k}\right)+\sum _{a=1}^{q-1}s^{a}=(-1)(-1)+(-1)=0.\)
			By Lam and Leung~\cite[Example~2.5]{LamLeung2000}, $V$ is a minimal vanishing sum: no nonempty proper subsum of $V$ vanishes. Since any reverse $\ell$-refinement ($\ell \in {r, p, q}$) among the singletons would contribute a proper $\ell$-cycle $\zeta^d \sum_{j=0}^{\ell-1} (\zeta^{rpq/\ell})^j = 0$ to $V$, no such refinement exists.
		\end{itemize}
		Thus $C_{r,p,q}$ is an irreducible ECS.
	\end{proof}
\begin{example} Setting $(r,p,q)=(2,3,5)$ in Theorem~A yields \[ \begin{aligned} C_{2,3,5} =\{&2(6),4(6),1(10),3(10),7(10),9(10),0(15),\\ &5(30),6(30),12(30),18(30),24(30),25(30)\}, \end{aligned} \] which is the classical irreducible ECS of length $30$ appearing in \cite{Schnabel2015}. \end{example} 
\begin{proposition} The system \[ \begin{aligned} \mathcal B =\{&1(6),5(6),4(10),6(10),0(15),3(15),12(15),\\ &2(30),8(30),9(30),10(30),20(30),21(30),22(30),28(30)\} \end{aligned} \] is an irreducible exact covering system of length $30$. Moreover, every modulus occurring in $\mathcal B$ has multiplicity at least two. 
\end{proposition} 
\begin{proof} 
		A direct inspection of the residues modulo $30$ shows that $\mathcal B$ partitions $[0,29]$. Since \[ \frac{2}{6}+\frac{2}{10}+\frac{3}{15}+\frac{8}{30}=1, \] it follows that $\mathcal B$ is an exact covering system of length $30$. To verify irreducibility, we rule out reverse prime refinements at every modulus:
		\begin{itemize}
			\item For modulus $6$, the residues $\{1, 5\}$ are incongruent modulo $3$ (ruling out a $2$-split) and fewer than $3$ classes appear (ruling out a $3$-split).
			\item For modulus $10$, the residues $\{4, 6\}$ are incongruent modulo $5$ and fewer than $5$ classes appear.
			\item For modulus $15$, the residues $\{0, 3, 12\}$ are pairwise incongruent modulo $5$ ($0, 3, 2 \pmod{5}$) and fewer than $5$ classes appear.
			\item For modulus $30$, reducing $R = \{2, 8, 9, 10, 20, 21, 22, 28\}$ modulo $15$, $10$, and $6$ yields $\{2, 8, 9, 10, 5, 6, 7, 13\}$, $\{2, 8, 9, 0, 0, 1, 2, 8\}$, and $\{2, 2, 3, 4, 2, 3, 4, 4\}$, respectively. No residue occurs twice modulo $15$, three times modulo $10$, or five times modulo $6$.
		\end{itemize}
		Hence no reverse prime refinement exists, and $\mathcal{B}$ is irreducible.
\end{proof} 
\begin{remark} Every ECS in the family $C_{r,p,q}$ contains a unique class of modulus $pq$, whereas every modulus occurring in $\mathcal B$ has multiplicity at least two. Therefore $\mathcal B$ does not belong to the family of Theorem~A. To the best of our knowledge, $\mathcal B$ is the first explicit irreducible ECS whose irreducibility is not explained by the Lam--Leung minimal vanishing sum construction. This suggests that the class of irreducible ECSs is considerably richer than the family provided by Theorem~A. 
\end{remark}
 \begin{remark} Additional infinite families of irreducible ECSs arise by taking any prime $s$, choosing distinct primes $r,p,q$ (one of which may equal $s$), and replacing the class $0(s)$ in the trivial partition \[ \{\,k(s):0\le k\le s-1\,\} \] by the expansion $sC_{r,p,q}$. The resulting ECS has length $srpq$ and remains irreducible.
  \end{remark}
	\section{Classification of Realizable Moduli Multisets}\label{sec:classification}
	We now formalize candidate multisets and establish Theorem~B.
	\begin{definition}\label{def:candidates}
		For every integer $n \geq 1$, let $\mathcal{C}(n)$ denote the set of all finite multisets $M = \{n_1, \dots, n_k\}$ of positive integers satisfying:
		\begin{enumerate}
			\item $\lcm(n_1, \dots, n_k) = n$;
			\item $D(M) = \sum_{i=1}^k \frac{1}{n_i} = 1$;
			\item $\gcd(n_i, n_j) > 1$ for all $1 \leq i < j \leq k$;
			\item If $n > 1$ and $p_1 < \cdots < p_m$ are the distinct prime divisors of $\max(M)$, then the multiplicity $t$ of $\max(M)$ in $M$ satisfies $t = \sum_{j=1}^m t_j p_j$ for some integers $t_j \geq 0$.
		\end{enumerate}
		Let $\mathcal{R}(n)\subseteq \mathcal{C}(n)$ denote the subset of multisets $M \in \mathcal{C}(n)$ realized by at least one ECS.
	\end{definition}

	We begin with prime-power lengths, where $\Par_p(m)$ denotes the number of partitions of $m$ into nonnegative powers of $p$.
\begin{proposition}\label{prop:primepower} For every prime $p$ and every integer $k\ge1$, \[ \mathcal R(p^k)=\mathcal C(p^k), \qquad |\mathcal R(p^k)|=\Par_p(p^k-1). \] \end{proposition}
	\begin{proof}
		Write a multiset $M$ of divisors of $p^k$ as $M = \{p^{(c_1)}, (p^2)^{(c_2)}, \dots, (p^k)^{(c_k)}\}$ with $c_j \geq 0$. The condition $\lcm(M) = p^k$ requires $c_k \geq 1$, and $D(M) = 1$ multiplies through by $p^k$ to give $\sum_{j=1}^k c_j p^{k-j} = p^k$, which automatically forces $p \mid c_k$ (hence $c_k \geq p$, satisfying~\eqref{eq:chen_porubsky}). Subtracting $1$ from both sides yields
		\(\sum _{j=1}^{k-1}c_{j}p^{k-j}+(c_{k}-1)p^{0}=p^{k}-1,\)
		establishing a bijection between $\mathcal{C}(p^k)$ and partitions of $p^k - 1$ into powers of $p$. Thus $|\mathcal{C}(p^k)| = \Par_p(p^k - 1)$.
		To realize $M \in \mathcal{C}(p^k)$, start from $\{0(1)\}$ and build a rooted $p$-ary refinement tree. Let $b_j$ denote the number of available classes of modulus $p^j$ before selecting $c_j$ of them as terminal leaves. Then $b_1 = p$ and $b_{j+1} = p(b_j - c_j)$, so by induction
		\(b_{j}-c_{j}=p^{j}\left(1-\sum _{i=1}^{j}\frac{c_{i}}{p^{i}}\right)=p^{j}\sum _{i=j+1}^{k}\frac{c_{i}}{p^{i}}\ge 0\)
		for every $1 \leq j \leq k$, with $b_k = c_k$. Thus every $M \in \mathcal{C}(p^k)$ is realized by an NECS.
	\end{proof}
	Next, we establish monotonicity of non-realizability under multiplication and two auxiliary lemmas.
	\begin{proposition}[Monotonicity of Non-Realizability]\label{prop:monotonicity}
		Let $n > 1$ be an integer. If $\mathcal{R}(n) \neq \mathcal{C}(n)$, then $\mathcal{R}(m) \neq \mathcal{C}(m)$ for every integer multiple $m$ of $n$.
	\end{proposition}
	\begin{proof}
		Let $M = \{n_1, \dots, n_k\} \in \mathcal{C}(n) \setminus \mathcal{R}(n)$ with $n_1 \leq \cdots \leq n_k$. Since $n > 1$ and $\lcm(M) = n$, we have $n_k \geq 2$. Let $m = sn$ for an integer $s \geq 2$, and define
		\(M^{\prime }=\{s^{(s-1)}\}\cup \{sn_{1},\dots ,sn_{k}\}.\)
		We verify that $M' \in \mathcal{C}(m)$:
		\begin{enumerate}
			\item $\lcm(M') = \lcm(s, sn_1, \dots, sn_k) = s\lcm(n_1, \dots, n_k) = sn = m$.
			\item $D(M') = \frac{s-1}{s} + \frac{1}{s}\sum_{i=1}^k \frac{1}{n_i} = \frac{s-1}{s} + \frac{1}{s} = 1$.
			\item Every element of $M'$ is divisible by $s > 1$, so $\gcd(x, y) \geq s > 1$ for all $x, y \in M'$.
			\item Since $n_k \geq 2$, we have $sn_k > s$, so $\max(M') = sn_k$ has the same multiplicity $t$ as $n_k$ in $M$. Because every prime divisor of $n_k$ divides $sn_k$, any representation $t = \sum_{\ell \mid n_k} t_\ell \ell$ ($t_\ell \in \mathbb{Z}{\geq 0}$) remains valid for the prime divisors of $sn_k$.
		\end{enumerate}
		Suppose for contradiction that $M' \in \mathcal{R}(m)$, realized by an ECS $\mathcal{A}'$. The $s-1$ classes of modulus $s$ in $\mathcal{A}'$ are pairwise disjoint and thus occupy $s-1$ distinct residue classes modulo $s$. Since every remaining class in $\mathcal{A}'$ has modulus $sn_i$ divisible by $s$, all remaining classes are contained in the unique uncovered residue class $r_0(s)$, where $0 \leq r_0 < s$. Writing each remaining class as $(r_0 + sa_i)(sn_i)$ with $0 \leq a_i < n_i$, exact coverage of $r_0(s)$ implies that $\{a_i(n_i)\}_{i=1}^k$ is an ECS realizing $M$, contradicting $M \notin \mathcal{R}(n)$. Thus $M' \in \mathcal{C}(m) \setminus \mathcal{R}(m)$.
	\end{proof}
	\begin{lemma}[Prime-Factorization Dichotomy]\label{lem:dichotomy}
		Every integer $n > 1$ either has the form $p^k$, $pq$, or $p^2 q$ for distinct primes $p \neq q$ and an integer $k \geq 1$, or is divisible by $pqr$, $p^3 q$, or $p^2 q^2$ for some pairwise distinct primes $p, q, r$.
	\end{lemma}
	\begin{proof}
		Let $n = \prod_{i=1}^{\omega(n)} p_i^{\alpha_i}$ be the prime factorization of $n > 1$, with $\alpha_i \geq 1$.
		\begin{itemize}
			\item If $\omega(n) \geq 3$, then $p_1 p_2 p_3 \mid n$.
			\item If $\omega(n) = 1$, then $n = p_1^{\alpha_1}$ with $\alpha_1 \geq 1$.
			\item If $\omega(n) = 2$, then $n = p_1^{\alpha_1} p_2^{\alpha_2}$. If $\alpha_1 + \alpha_2 \leq 3$, then $\{\alpha_1, \alpha_2\} \in \{\{1,1\}, \{2,1\}\}$, so $n$ has the form $pq$ or $p^2 q$. If $\alpha_1 + \alpha_2 \geq 4$, then either $\max(\alpha_1, \alpha_2) \geq 3$ (so $p^3 q \mid n$ for $\{p, q\} = \{p_1, p_2\}$) or $\alpha_1, \alpha_2 \geq 2$ (so $p_1^2 p_2^2 \mid n$). \qedhere
		\end{itemize}
	\end{proof}
	\begin{lemma}[Modular Packing Bound]\label{lem:pigeonhole}
		Let $p$ and $q$ be distinct primes. Among arithmetic progressions of modulus $pq$, at most $q$ pairwise disjoint progressions can be contained in a single residue class modulo $p$.
	\end{lemma}
	\begin{proof}
		Fix a residue class $r_0(p)$ with $0 \leq r_0 < p$. Any arithmetic progression $a(pq)$ contained in $r_0(p)$ satisfies $a \equiv r_0 \pmod{p}$, so its residue modulo $pq$ belongs to $\{r_0 + jp : 0 \leq j \leq q - 1\}$, which has cardinality $q$. By the Pigeonhole Principle, any collection of more than $q$ progressions of modulus $pq$ contained in $r_0(p)$ contains two progressions with the same residue modulo $pq$, which coincide.
	\end{proof}
	\begin{proof}[Proof of Theorem~B]
		We first verify that $\mathcal{R}(n) = \mathcal{C}(n)$ whenever $n$ has the form $p^k$, $pq$, or $p^2 q$ for distinct primes $p \neq q$ and $k \geq 1$, and then exhibit non-realizable candidates at each minimal remaining prime signature.
		\medskip
		\noindent\emph{Step 1: Realizability for $p^k$, $pq$, and $p^2 q$.}
		The equality $\mathcal{R}(p^k) = \mathcal{C}(p^k)$ is Proposition~\ref{prop:primepower}.
		For $n = pq$, let $M \in \mathcal{C}(pq)$. Since $\lcm(M) = pq > 1$, we have $1 \notin M$, so $M = \{p^{(a)}, q^{(b)}, (pq)^{(c)}\}$ with $c \geq 1$. Since $\gcd(p, q) = 1$, the pairwise non-coprime condition forces $a = 0$ or $b = 0$. Without loss of generality, assume $b = 0$. Multiplying $D(M) = 1$ by $pq$ gives $aq + c = pq$, whence $c = (p - a)q$ with $0 \leq a < p$. Starting from $\{0(p), 1(p), \dots, (p-1)(p)\}$ and applying a $q$-refinement to $p - a$ of its classes realizes $M$.
		For $n=p^2q$, let $M\in \mathcal C(p^2q)$. Since $\lcm(M)=p^2q>1$, we have $1\notin M$, so \[ M\subseteq \{p,p^2,q,pq,p^2q\}. \] Because $\gcd(q,p)=\gcd(q,p^2)=1$, the pairwise non-coprime condition implies that $q$ cannot occur together with either $p$ or $p^2$. Hence exactly one of the following two cases occurs.
		
		 \medskip \noindent\emph{Case A:} \[ M\subseteq\{q,pq,p^2q\}. \] Write \[ M=\{q^{(a)},(pq)^{(b)},(p^2q)^{(d)}\}. \] Dividing all moduli by $q$ yields \[ N=\{1^{(a)},p^{(b)},(p^2)^{(d)}\}. \] Since $D(M)=1$, $D(N)=q$. Thus \[ a+\frac{b}{p}+\frac{d}{p^2}=q. \] Multiplying by $p^2$ gives \[ bp+d=(q-a)p^2. \] Hence $p\mid d$, say $d=up$. Then \[ b+u=(q-a)p. \] Consequently $N$ decomposes into $a$ copies of $\{1\}$ and $q-a$ multisets \[ N_j=\{p^{(b_j)},(p^2)^{(u_jp)}\}, \qquad b_j+u_j=p. \] Each $N_j$ belongs to $\mathcal C(p^2)$ and is realizable by Proposition~\ref{prop:primepower}. If $\mathcal A_j$ realizes $N_j$, then \[ \mathcal A = \bigcup_{j=0}^{q-1}(q\mathcal A_j+j) \] is an ECS realizing $M$.
		 
		  \medskip \noindent\emph{Case B:} \[ M\subseteq\{p,p^2,pq,p^2q\}. \] Write \[ M=\{p^{(a)},(p^2)^{(b)},(pq)^{(c)},(p^2q)^{(d)}\}. \] By~\eqref{eq:chen_porubsky}  applied to the largest modulus $p^2q$, there exist $t_1,t_2\ge0$ such that \[ d=t_1p+t_2q . \] Define \[ M' = \{p^{(a)},(p^2)^{(b+t_2)},(pq)^{(c+t_1)}\}. \] Then \[ D(M')=D(M)=1. \] Dividing all moduli of $M'$ by $p$ gives \[ N' = \{1^{(a)},p^{(b+t_2)},q^{(c+t_1)}\}, \] with $D(N')=p$. Therefore \[ a+\frac{b+t_2}{p}+\frac{c+t_1}{q}=p. \] Since $\gcd(p,q)=1$, \[ p\mid (b+t_2), \qquad q\mid (c+t_1). \] Write \[ b+t_2=v_1p, \qquad c+t_1=v_2q. \] Then \[ a+v_1+v_2=p. \] Hence $N'$ decomposes into $a$ copies of $\{1\}$, $v_1$ copies of $\{p^{(p)}\}$, and $v_2$ copies of $\{q^{(q)}\}$. Each of these multisets is realizable by Proposition~\ref{prop:primepower}, so $N'$ is realizable. Expanding by $p$ shows that $M'$ is realizable. Finally, perform $t_2$ $q$-refinements on classes of modulus $p^2$ and $t_1$ $p$-refinements on classes of modulus $pq$. This replaces $t_2$ classes of modulus $p^2$ by $t_2q$ classes of modulus $p^2q$ and $t_1$ classes of modulus $pq$ by $t_1p$ classes of modulus $p^2q$. The resulting ECS realizes $M$.
		  
		\medskip
		\noindent\emph{Step 2: Non-realizability for $p^3 q$, $p^2 q^2$, and $pqr$.}
		We construct an explicit multiset $M \in \mathcal{C}(L) \setminus \mathcal{R}(L)$ for each of the three minimal non-realizable forms:
		\begin{itemize}
			\item \emph{$L = p^3 q$ ($p \neq q$):} Define
		\[ M= \begin{cases} \{\,p^{(p-1)},\,(p^2)^{(p-1)},\,pq,\,(p^3q)^{(p(q-p))}\,\}, & \text{if } p<q,\\[4pt] \{\,p^{(p-1)},\,p^3,\,(pq)^{(q-1)},\,(p^3q)^{(p^2-q)}\,\}, & \text{if } p>q. \end{cases} \]
			Let $c$ denote the multiplicity of $p^3 q$. If $p < q$, then $c = p(q - p) \geq p \geq 2$ is a positive multiple of $p$; if $p > q$, then writing $p = q + d$ with $d \geq 1$ yields $c = p^2 - q = p(q + d) - q = (d - 1)p + (p - 1)q$, where $d - 1 \geq 0$ and $p - 1 \geq 2$. In both cases, $c$ satisfies~\eqref{eq:chen_porubsky}, $\lcm(M) = p^3 q$, every element is divisible by $p$, and $D(M) = 1$, so $M \in \mathcal{C}(p^3 q)$. In any ECS $\mathcal{A}$ realizing $M$, the $p - 1$ classes of modulus $p$ occupy $p - 1$ residue classes modulo $p$, forcing a class $a(p^k)$ ($k \in \{2, 3\}$) and a class $b(pq)$ into the single remaining residue class modulo $p$. Then $p \mid (a - b)$, contradicting Lemma~\ref{lem:intersection} since $\gcd(p^k, pq) = p$. Thus $M \notin \mathcal{R}(p^3 q)$.
			\item \emph{$L = p^2 q^2$ ($p < q$):} Define
			\(M=\{p^{(p-1)},\;(p^{2})^{(p-1)},\;pq,\;(p^{2}q^{2})^{(q(q-p))}\}.\)
			Since $q > p \geq 2$, $c = q(q - p) \geq q > p$ is a positive multiple of $q$, $\lcm(M) = p^2 q^2$, every element is divisible by $p$, and $D(M) = 1$, so $M \in \mathcal{C}(p^2 q^2)$. By the same obstruction as above, the $p - 1$ classes of modulus $p$ force a class $a(p^2)$ and a class $b(pq)$ into the same residue class modulo $p$, contradicting Lemma~\ref{lem:intersection}. Thus $M \notin \mathcal{R}(p^2 q^2)$.
			\item \emph{$L = pqr$ ($p < q < r$):} Since $\gcd(p,q)=1$, multiplication by $q$ induces a permutation of the residue classes modulo~$p$. Therefore there exists a unique $b\in\{1,\ldots,p-1\}$ such that \[ bq\equiv (q-1)r \pmod p. \] Define \[ M=\{(pq)^{((p-1)q+1)},\;(pr)^{(b)},\;(pqr)^{(c)}\}, \qquad c=(q-1)r-bq. \] Since $q\ge p+1$ and $r>q$, we have \[ (q-1)r\ge pr>pq>(p-1)q\ge bq, \] so $c>0$. By the choice of $b$, \[ c=(q-1)r-bq\equiv0\pmod p, \] so $c$ is a positive multiple of $p$ and satisfies~\eqref{eq:chen_porubsky} for $pqr$. Moreover, $\lcm(M)=pqr$, every element is divisible by $p$, and \[ D(M) = \frac{(p-1)q+1}{pq} +\frac{b}{pr} +\frac{(q-1)r-bq}{pqr} =1, \] so $M\in\mathcal{C}(pqr)$. If $\mathcal{A}$ realized $M$, then since $\mathcal{A}$ contains $(p-1)q+1$ pairwise disjoint classes of modulus $pq$, Lemma~\ref{lem:pigeonhole} and the Pigeonhole Principle imply that every residue class modulo $p$ contains at least one class of modulus $pq$. Hence a class $a(pr)\in\mathcal{A}$ must share a residue class modulo $p$ with some class $d(pq)\in\mathcal{A}$, so $p\mid(a-d)$. Since $\gcd(pr,pq)=p$, Lemma~\ref{lem:intersection} implies that \[ a(pr)\cap d(pq)\neq\varnothing, \] contradicting the assumption that $\mathcal{A}$ is an exact covering system. Thus $M\notin\mathcal{R}(pqr)$.
			 \end{itemize}
		\medskip
		\noindent\emph{Step 3: Conclusion.}
		By Lemma~\ref{lem:dichotomy}, every integer $n > 1$ not of the form $p^k$, $pq$, or $p^2 q$ is a multiple of $p^3 q$, $p^2 q^2$, or $pqr$ for some pairwise distinct primes $p, q, r$. Applying Step~2 and Proposition~\ref{prop:monotonicity} shows that $\mathcal{R}(n) \neq \mathcal{C}(n)$ for all such $n$, completing the proof.
	\end{proof}
	
\end{document}